\documentclass[10pt]{amsart} %
\usepackage{fullpage}
\usepackage{graphicx}
\usepackage{graphics}
\usepackage{amsmath}
\usepackage{amsfonts}
\usepackage{amssymb}
\usepackage{amsthm}
\usepackage{psfrag}
\usepackage{amsmath,amssymb,color}
\usepackage{amsfonts,amstext}
\usepackage{enumerate}

\usepackage[normalem]{ulem} 
\usepackage{soul}

\newtheorem{theorem}{Theorem}

\newtheorem{conjecture}{Conjecture}

\title{Strong edge-colorings for $k$-degenerate graphs}
\date{\today}

\author{Gexin Yu}
\thanks{The author's research was supported in part by an NSA grant. }

\address{Department of Mathematics, College of William and Mary, Williamsburg, VA, 23185. }\email{gyu@wm.edu}

\begin{document}
\maketitle

\begin{abstract}
We prove that the strong chromatic index for each $k$-degenerate graph with maximum degree $\Delta$ is at most $(4k-2)\Delta-k(2k-1)+1$. 
\end{abstract}


A {\em strong edge-coloring} of a graph $G$ is an edge-coloring so that no edge can be adjacent  to two edges with the same color. So in a strong edge-coloring, every color class gives an induced matching. The strong chromatic index $\chi_s'(G)$ is the minimum number of colors needed to color $E(G)$ strongly.  This notion was introduced by Fouquet and Jolivet (1983, \cite{FJ83}).   Erd\H{o}s and Ne\v{s}et\v{r}il during a seminar in Prague in 1985 proposed some open problems, one of which is the following

\begin{conjecture}[Erd\H{o}s and Ne\v{s}et\v{r}il, 1985]
If $G$ is a simple graph with maximum degree $\Delta$, then $\chi_s'(G)\le 5\Delta^2/4$ if $\Delta$ is even, and $\chi_s'(G)\le (5\Delta^2-2\Delta+1)/4$ if $\Delta$ is odd.
\end{conjecture}

This conjecture is true for $\Delta\le 3$ (\cite{A92, HHT93}). Cranston \cite{C06} showed that $\chi_s'(G)\le 22$ for $\Delta=4$.   Chung, Gy\'arf\'as, Trotter, and Tuza (1990, \cite{CGTT90}) showed that the upper bounds are exactly the numbers of edges in $2K_2$-free graphs.  Molloy and Reed \cite{MR97} proved that graphs with sufficient large maximum degree $\Delta$ has strong chromatic index at most $1.998\Delta^2$. For more results see \cite{SSTM} (Chapter 6, problem 17). 


A graph is {\em $k$-degenerate} if every subgraph has minimum degree at most $k$.  Chang and Narayanan (2012, \cite{CN12}) recently proved that a $2$-degenerate graph with maximum degree $\Delta$ has strong chromatic index at most $10\Delta-10$.  Luo and the author in \cite{LY12} improved the upper bound to $8\Delta-4$.  

In~\cite{CN12}, the following conjecture was made

\begin{conjecture}[Chang and Narayanan, \cite{CN12}]
There exists an absolute constant $c$ such that for any $k$-degenerate graphs $G$ with maximum degree $\Delta$, $\chi_s'(G)\le ck^2\Delta$.  Furthermore, the $k^2$ may be replaced by $k$. 
\end{conjecture}

In this paper, we prove a stronger form of the conjecture.  Unlike the priming processes in\cite{CN12, LY12}, we find a special ordering of the edges and by using a greedy coloring obtain the following result. 

\begin{theorem}
 The strong chromatic index for each $k$-degenerate graph with maximum degree $\Delta$ is at most $(4k-2)\Delta-k(2k-1)+1$. 
 \end{theorem}
 
 Thus,  $2$-degenerate graphs have strong chromatic index at most $6\Delta-5$.

\begin{proof}
By definition of $k$-degenerate graphs, after the removal of all vertices of degree at most $k$, the remaining graph has no edges or has new vertices of degree at most $k$, thus we have the following simple fact on $k$-degenerate graphs (see also \cite{CN12}). 

\medskip

{\em Let $G$ be a $k$-degenerate graph.  Then there exists $u\in V(G)$ 
 so that $u$ is adjacent to at most $k$ vertices of degree more than $k$.  Moreover, if $\Delta(G)>k$, then the vertex $u$ can be selected with degree more than $k$.}

\medskip

We call a vertex $u$ 
a {\em special vertex} if $u$ is adjacent to at most $k$ vertices of degree more than $k$.  An edge is a {\em special edge} if it is incident to a special vertex and a vertex with degree at most $k$.  The above fact implies that every $k$-degenerate graph has a special edge, and if $\Delta\le k$, then every vertex and every edge are special. 

We order the edges of $G$ as follows.   First we find in $G$ a special edge,  put it at the beginning of the list, and then remove it from $G$.   Repeat the above step in the remaining graph.  When the process ends, we have an ordered list of the edges in $G$, say $e_1, e_2, \ldots, e_m$, where $m=|E(G)|$.  So $e_m$ is the special edge we first chose and placed in the list. 

Let $G_i$ be the graph induced by the first $i$ edges in the list, $i=1,2,\ldots, m$.  Then $e_i$ is a special edge in $G_i$. 
We now count the edges of $G_i$ within distance one to $e_i$ in $G$.  We may call the edges in $G_i$ blue edges and the edges in $G-G_i$ yellow edges.  Let $u_i,v_i$ be the endpoints of $e_i$ with $u_i$ being a special vertex in $G_i$.  


We first count the blue edges incident to $u_i$ and its neighbors.  The vertex $u_i$ has three kinds of neighbors:   the neighbors in $X_1$ sharing blue edges with $u_i$ and having degree more than $k$,  the neighbors in $X_2$ sharing blue edges with $u_i$ and having degree at most $k$ (thus $v_i\in X_2$), and the neighbors in $X_3$ sharing yellow edges with $u_i$.  By definition, $|X_1|\le k$, so  at most $|X_1|\Delta+k(|X_2|-1)$ blue edges  are incident to $X_1\cup (X_2-\{v_i\})$.  For each vertex $u$ in $X_3$, $uu_i$ is a yellow edge in $G_i$ but will be a special edge in $G_j$ for some $j>i$.  So either $u$ or $u_i$ has degree at most $k$ in $G_j$ (thus also in $G_i$), and if $u_i$ has degree at least $k$ in $G_m$ for some $m$, then all yellow edges incident to $u_i$ in $G_m$ should have degree at most $k-1$ in $G_m$, in order for the yellow edges to be special later.  Then among vertices in $X_3$, at most $x=\max\{0,k-|X_1|-|X_2|\}$ vertices have degree more than $k$ in $G_i$, and all other vertices have degree at most $k-1$ in $G_i$.  Therefore at most $x\Delta+(|X_3|-x)(k-1)$ blue edges are incident to $X_3$.    Note that $d(u_i)\le \Delta, |X_2|\le \Delta$ and $|X_1|+x\le k$, then at most 
$$|X_1|\Delta+k(|X_2|-1)+x\Delta+(|X_3|-x)(k-1)=(|X_1|+x)\Delta+(k-1)(d(u_i)-|X_1|-x-1)+|X_2|-1\le 2k\Delta-k^2$$ 
blue edges are within distance one to $e_i$ from $u_i$ side (not including the edges incident to $v_i$).

We also count the blue edges incident to $v_i$ and its neighbors. Similarly, $v_i$ has two kinds of neighbors:  the neighbors in $Y_1$ sharing blue edges with $v_i$, and the neighbors in $Y_2$ sharing yellows edges with $v_i$.  From the fact that $e_i$ is a special edge, $|Y_1|\le k$, so at most $(|Y_1|-1)\Delta$ blue edges are incident to $Y_1-\{u_i\}$.  For each vertex $v$ in $Y_2$, $vv_i$ is a yellow edge in $G_i$ but will be a special edge in $G_s$ for some $s>i$.  Similar to above,  at most $k-|Y_1|$ vertices in $Y_2$ have degree more than $k$ in $G_i$, and all other vertices in $Y_2$ have degree at most $k-1$ in $G_i$.  So at most $(k-|Y_1|)(\Delta-1)+(|Y_2|-(k-|Y_1|))(k-1)$ blued edges are incident to $Y_2$. In total, at most 
$$(|Y_1|-1)\Delta+(k-|Y_1|)(\Delta-1)+(|Y_2|-(k-|Y_1|))(k-1)\le (2k-2)\Delta-k(k-1)$$


So in $G_i$, the number of blue edges within distance one to $e_i$ is at most 
$$2k\Delta-k^2+(2k-2)\Delta-k(k-1)\le (4k-2)\Delta-k(2k-1)$$

Now color the edges in the list one by one greedily.  For each $i$, when it is the turn to color $e_i$, only the edges in $G_i$ (the blue edges) have been colored.  Since there are at least $(4k-2)\Delta-k(2k-1)+1$ colors,  we are able to color the edges so that edges within distance one get different colors.  
\end{proof}

We shall note that the above result is not only true for simple graphs, but also for multigraphs.  

\section*{Acknowledgement}

The author would like to thank Rong Luo and Zixia Song for their encouragements and discussions.

\end{document}